\input amstex 
\documentstyle{amsppt}
\input bull-ppt
\keyedby{bull424e/kmt}

\define\Fhat{{\widehat F}}
\define\Ftilde{{\widetilde F}}
\define\SC{Sierpinski carpet }
\define\R{{\Bbb R}}

\redefine\P{{\Bbb P}}
\define\Z{{\Bbb Z}}
\redefine\S{{\Cal S}}

\define\bz{{\bold z_0}}

\define\q {\quad} \define\qq {\qquad}

\topmatter
\cvol{29}
\cvolyear{1993}
\cmonth{October}
\cyear{1993}
\cvolno{2}
\cpgs{208-212}
\ratitle
\title Coupling and  Harnack inequalities  \\ for 
Sierpinski carpets \endtitle
\author Martin T. Barlow  and Richard F. Bass \endauthor
\shortauthor{M. T. Barlow and R. F. Bass}
\address (Martin T. Barlow) Department of Mathematics, 
University of British Columbia, 
Vancouver,  British Columbia V6T 1Z2, Canada\endaddress
\ml barlow\@unixg.ubc.ca\endml  
\address (Richard F. Bass) Department of Mathematics, 
University of Washington, Seattle,
Washington 98195\endaddress
\ml bass\@math.washington.edu \endml
\date November 9, 1992\enddate
\subjclass Primary 60B99;  Secondary 60J35\endsubjclass
\thanks Research partially supported by NSF grant DMS 
91--00244 \endthanks
\keywords Harnack inequality, Sierpinski carpets, 
self--similar, fractals,
Brownian motion, heat equation, transition densities, 
Poincar\'e
inequality, Sobolev inequality, spectral dimension, 
electrical
resistance\endkeywords
\abstract Uniform Harnack inequalities for harmonic 
functions on the pre- and
graphical  Sierpinski carpets are proved using a 
probabilistic coupling
argument. Various  results follow from this, including the 
construction of
Brownian motion on  Sierpinski carpets embedded in $\R^d$, 
$d\geq 3$, estimates
on the fundamental  solution of the heat equation, and 
Sobolev and Poincar\'e
inequalities. \endabstract

\endtopmatter

\document

The Sierpinski carpets (SCs) we will study are 
generalizations
of the Cantor set. Let $F_0 = [0,1]^d$ be the unit cube in 
$\R^d$, $d\geq 2$, centered at  $\bz=(1/2,\ldots,1/2)$. 
Let $k,a$ be 
integers with $1\le a <k$ and $a+k$ even. 
Divide $F_0$ into $k^d$ equal subcubes, remove a central 
block of
$a^d$ subcubes, and let $F_1$ be what remains: thus $F_1=F_0
-((k-a)/2k,(k+a)/2k)^d$.
Now repeat this operation on each of the $k^d-a^d$ 
remaining subcubes to
obtain
$F_2$. Iterating, we obtain a decreasing sequence of 
closed sets
$F_n$; then $F=
\bigcap_{n=0}^\infty F_n$ is a \SC and has Hausdorff 
dimension
$d_f=d_f(F)=
\log(k^d-a^d)/\log(k)$.
(When $d=2$, $k=3$, and $a=1$, we get the
usual Sierpinski carpet.) Let $\Fhat_n=k^nF_n \subset 
[0,\infty)^d$, and define 
the {\it pre-Sierpinski carpet} by 
$\Fhat=\bigcup_{n=1}^\infty \Fhat_n$
(see \cite{10}).
The {\it graphical Sierpinski carpet} is the graph $G=(V,E)$
with vertex set $V=(\bz+\Z^d) \cap \Fhat$ and edge set 
$E=\{\, \{x,y\} \in V: \vert x-y \vert =1 \}$.

Thus $\hbox{int}(\Fhat)$ is a domain in $\R^d$ with a 
large-scale structure
which mimics the small-scale structure of $F$. We are 
interested in the
behavior of solutions of the Laplace and heat equations on 
$F$, $\Fhat$,
and $G$. One reason for this is applications to ``transport
phenomena'' in disordered media (see \cite{6}); another is 
the new type of
behavior of the heat kernel on these spaces. Let $W$ be 
Brownian
motion on $\Fhat$ with normal reflection on $\partial 
\Fhat$, and 
let $q(t,x,y)$ be the\ transition 
density of $W$, so that $q$ solves
the heat equation on $\Fhat$ with Neumann boundary 
conditions\ on
$\partial \Fhat$.

\proclaim {Theorem 1} There exist $c_1, \ldots, c_6 \in 
(0,\infty)$
and $d_s=d_s(F) \in (1,d_f)$ such that
if $x,y\in \Fhat$, $t\in (1,\infty)$, $|x-y| \le t$, then
$$\aligned
&c_1 t^{-d_s/2} 
\exp\bigg(-c_2 
\Big(\frac{|x-y|^{d_w}}{t}\Big)^{1/(d_w-1)}\bigg) \\
&\qquad\leq q(t,x,y) \leq
c_3 t^{-d_s/2} \exp\bigg(-c_4\Big(\frac{|x-y|^{d_w}}{t}%
\Big)^{1/(d_w-1)}\bigg),
 \endaligned\tag 1$$
where $d_w=2d_f/d_s$;
while if $x,y\in \Fhat$, $t\in (1,\infty)$, $|x-y| > t$, 
then
$$
\exp\left(-c_5 \frac{|x-y|^2}{t}\right) \le q(t,x,y) \le 
\exp\left(-c_6 
\frac{|x-y|^2}{t}\right).
\tag2
$$
\endproclaim

The index $d_s$ is called the {\it spectral dimension} of 
$F$ and
turns out to be much more significant than the Hausdorff 
dimension $d_f$
as far as analytic properties of these spaces are concerned.
Since $d_s<d$,  (1) confirms the physical intuition that 
the presence
of increasingly large 
reflecting barriers causes heat to dissipate to infinity 
more slowly.
It seems unlikely that there is any simple relationship 
between $d_s, k, a$,
and $d$.

While there is a well-developed
approach to the heat equation  using analytic tools such as
Sobolev or log-Sobolev inequalities (see \cite{7}), these
methods do not appear to give the best-possible results on 
spaces such 
as $\Fhat$---compare the upper bound on $q(t,x,y)$ given in
Theorem 1 with the results of \cite{10}. 

The proof of Theorem 1 rests on the following Harnack 
inequality.
Let $D \subset \R^d$ be open: we will say that $h$ is {\it 
harmonic on }
$D \cap \Fhat$ if (i) $\Delta h=0$ in 
$\hbox{int}(D \cap \Fhat)$ and (ii) $h$
has $0$ normal derivative a.e.\  on $D\cap \partial \Fhat$.
Equivalently, $h$ is harmonic with respect to $W(t \wedge 
T_D)$, where
$T_D = \inf \{t: W(t)\notin D \}$.
Let $D_n=(-1,k^n)^d$. 

\proclaim {Theorem 2} There exists $c_1\in (0,\infty)$ 
\RM(depending only 
on $d,k,a)$, such that
if $h$ is positive harmonic in $D_n \cap \Fhat$ and 
$x,y\in D_{n-1} \cap \Fhat$, then $h(x)/h(y)\leq c_1$.
\endproclaim
\rem{Remarks} 1.
Note that $c_1$ is independent of $n$; otherwise the
result is trivial. 
\par 2. The case $d=2$ was proved in
\cite{1}; the proof there relies on the fact that a closed 
curve
in the plane separates the plane into two pieces. Just as 
in the
case of elliptic operators, the results for two dimensions 
are
considerably easier to prove. The result of \cite{1} was 
extended
in \cite{8} to SCs with $d_s(F) <2$.
\par  3. Using the symmetry of $\Fhat$, Theorem 2 extends 
to other
domains in $\Fhat$.
\par  4. Theorems 1 and 2 actually hold for a much wider 
class of 
SCs, those
satisfying a higher-dimensional generalization of (2.1) of 
\cite{4}.
\par 5. A similar result holds for the graphical
Sierpinski carpet $G$.
\par 6. 
Most existing proofs of Harnack inequalities for selfadjoint
operators depend on Sobolev inequalities, which in turn 
depend on the
underlying geometry of the space. Here the appropriate 
Sobolev inequality
involves the spectral dimension $d_s(F)$;
however, no geometric definition
of $d_s$ is known.
Thus we were led to abandon analytic approaches
in favor of the probabilistic coupling argument described 
at the end of
this paper.
\endrem

We now describe some other consequences of Theorem 2.  
Let $\Ftilde\!=\!\bigcup_{n=0}^\infty k^nF$\!, the SC
extended to $[0,\infty)^d$, and\ write 
$\mu$ for Hausdorff $x^{d_f}$-measure
on $\Ftilde$.

\proclaim {Theorem 3} There exists a strong Markov process
$X_t$ with state space $\Ftilde$ such that $X$ has a 
strong Feller transition
semigroup $P_t$ which is $\mu$-symmetric,
$X_t$ has continuous
paths, $X_t$ is self-similar with respect to dilations
of size $k^n$, and the process $X$ is locally invariant 
with respect
to the local isometries of $\Ftilde$.
\endproclaim

Let $p(t,x,y)$ be the transition density of $X_t$ with 
respect
to $\mu$. Then $p(t,x,y)$ is the fundamental solution to the
heat equation on $\Ftilde$: $\partial u/\partial 
t=\Delta_\Ftilde u$,
where $\Delta_\Ftilde$ is the infinitesimal 
generator of $X_t$. Then we have

\proclaim {Theorem 4} There exist $c_1, c_2, c_3, c_4 \in 
(0,\infty)$
and $d_s=d_s(F) \in (1,d_f)$ such that
for all $x,y\in \Ftilde$, $t\in (0,\infty)$,
$$\eqalign{&c_1 t^{-d_s/2} 
\exp\bigg(-c_2 
\Big(\frac{|x-y|^{d_w}}{t}\Big)^{1/(d_w-1)}\bigg)\cr
&\qquad\leq p(t,x,y) \leq
c_3 t^{-d_s/2} \exp\bigg(-c_4\Big(\frac{|x-y|^{d_w}}{t}%
\Big)^{1/(d_w-1)}\bigg),
\cr}
$$
where $d_w=2d_f/d_s$. Moreover, $p(t,x,y)$ is $C^\infty$ 
in $t$, and $p(t,x,y)$ and all
its partial derivatives with respect to $t$ are jointly 
H\"older continuous
in $x$ and $y$.
\endproclaim
Many properties of the process $X$, such as
its transience or recurrence, the existence of local 
times, the existence of self-intersections, and the 
asymptotic frequency of
eigenvalues follow easily from Theorem 4. 
For example, note that $X$ is point recurrent if and only if
$d_s(F) <2$.

The next set of consequences include Sobolev inequalities, 
Poincar\'e
inequalities, and electrical resistance inequalities for 
$\Ftilde$,
$\Fhat$, and $G$---nine theorems in total.
Since the electrical resistance  inequalities are probably 
the least
well-known type, we give the one for $G$ as a 
representative sample. 
If $B$ is any subset of $G$, let $|B|$ denote the 
cardinality
of $B$.
Then $R(B)$, the resistance  from $B$ to infinity, is 
defined by
$$R(B)^{-1} = \inf \left\{ 
\sum_{ \{x,y\} \in E(G) } (f(x)-f(y))^2
: f\equiv 1 \hbox { on }
B, f(x)\to 0 \hbox{ as } |x|\to \infty \right\}.$$
The inverse of $R(B)$ is the conductance from $B$ to 
infinity and
equals the capacity of $B$.

\proclaim {Theorem 5} Suppose $d_s=d_s(F)>2$, and let 
$\zeta=d_s/(d_s-2)$. 
Then there exists $c_1$ such that if $A\subset G$,
$|A|\leq c_1 R(A)^{-\zeta}$.
\endproclaim

Theorem 5 follows fairly straightforwardly from Theorem 4 by
applying ideas of \cite{11} and \cite{12}. As Theorems 1, 
3, and 4 follow
from Theorem 2 by generalizations and modifications of 
methods of
\cite{1--4, 8, 9},
we discuss only Theorem 2.

Let $W_t$ be the Brownian motion on $\Fhat$ described 
above, and let
$$\tau(x,r)=\inf\{t:|W_t-x|\geq r\},\qquad
T(x,r)= \inf\{t:|W_t-x|\leq r\}.$$
The following lemma is proved in a similar fashion to 
Lemma 3.2 of 
\cite{1}.

\proclaim {Lemma 6} There exist $c_2>c_1>1$, $\delta>0$
independent of $r$ such that if $x,y\in \Fhat$ and
$|y-x|\leq c_1 r$, then
$$\P^y(T(x,r)<\tau(x,c_2r))>\delta. \tag3$$
\endproclaim

It is known (see
Theorem 3.9 of \cite{5}, for example) that the Harnack 
inequality
Theorem 2 follows from (3) and an oscillation inequality
of the following form. 

\proclaim {Lemma 7} There exists
$\rho <1$ such that if 
$n \ge 1$ and $h$ is positive harmonic on $D_n \cap \Fhat$,
then 
$$|h(x)-h(y)|\leq \rho \sup_{z\in D_n \cap \Fhat} |h(z)|,
\qq x,y\in D_{n-1}\cap \Fhat. \tag4$$
\endproclaim

To show (4), it suffices  to construct two
$\Fhat$-valued Brownian motions $W^x$ and $W^y$, starting 
from $x$ 
and $y$ respectively,
which couple (i.e., meet) with probability at least
$1-\rho$ before either exits $D_n$.

Fix $n$. Let $\S_m$ be the collection of cubes of side 
length $k^m$
with vertices in $k^m\Z^d$. 
Say that $x, y \in \Fhat$ are $m$-associated 
if there is an isometry of the cube in $\S_{m}$ containing
$x$ onto the cube in $\S_{m}$ containing $y$ that maps $x$ 
onto
$y$. Note that if two points are $m$-associated, then they 
will
also be $\ell$-associated for all $\ell \le m$. 

Suppose first that $x$ and $y$ are $m$-associated. 
We start a Brownian motion $W^x(t)$ on $\Fhat$ at $x$.
Let $U_0=0$, and $U_{i+1} =\inf\{t:|W^x(t)-W^x(U_i)| \geq 
k^{m}\}$.
The key step is to exploit the local symmetry of $\Fhat$ 
to construct, using suitable reflections, another Brownian 
motion $W^y(t)$ on
$\Fhat$, starting at $y$, such that
(a) $W^x(t)$ and $W^y(t)$ are
$m$-associated for all $t \ge 0$, and (b)
there exist $j$ and $c_1>0$ such that 
$$\P\Big(W^x(U_j(\omega))
\hbox { and } W^y(U_j(\omega)) 
\hbox { are $(m+1)$-associated }\Big) >c_1. \tag5$$ 

A renewal argument and then an induction show that if $x$ 
and $y$ are
0-associated, then $W^x$ and $W^y$ couple with probability 
$c_2>0$ 
before either process leaves $D_n$. Lemma 6 then follows 
easily.

\enddocument